\documentclass[10pt,twocolumn]{article}
\usepackage{amsmath,amssymb,amsfonts,mathtools}
\usepackage{bm}
\usepackage{booktabs}
\usepackage{tabularx}
\usepackage{array}
\usepackage{xcolor}
\usepackage{graphicx}
\usepackage{tikz}
\usepackage{microtype}
\usepackage{geometry}
\usepackage{titlesec}

\titleformat{\section}
  {\normalfont\large\bfseries}
  {\thesection}{0.5em}{}

\titleformat{\subsection}
  {\normalfont\normalsize\bfseries}
  {\thesubsection}{0.5em}{}
  
\usetikzlibrary{arrows.meta,positioning,calc,fit,shapes.geometric}

\newcolumntype{Y}{>{\raggedright\arraybackslash}X}

\newcommand{\data}{\mathcal{D}}
\newcommand{\loss}{\mathcal{L}}

\makeatletter
\renewcommand{\@seccntformat}[1]{\csname the#1\endcsname\hspace{0.35em}}
\makeatother

\begin{document}

\title{\bf From Horner's rule to neural networks}

\author{
Lorenzo Pareschi\footnote{Maxwell Institute and School of Mathematical Sciences, Heriot-Watt University, Edinburgh, UK and
Department of Mathematics and Computer Science, University of Ferrara, Italy.}
}

\date{}

\maketitle

\begin{abstract}
Horner's rule is usually taught as an efficient algorithm for evaluating a polynomial. Read as a recurrence, however, it is also a small computational graph: an input is propagated through intermediate states, while parameters enter locally. This note follows that elementary observation from polynomial evaluation to polynomial fitting. The equation $Va\simeq y$ then unifies least-squares approximation, interpolation and overparametrisation as different regimes of a fixed representation. 
Neural networks are introduced as the next step: the representation itself becomes parametrised and is learned from data. The purpose is pedagogical rather than encyclopaedic. It offers a numerical route into machine learning, using approximation, linear algebra and optimisation as familiar entry points.
\end{abstract}


\section{The polynomial as a model of learning}
\subsection{A classical algorithm as a computational graph}
There are algorithms so familiar that they almost disappear from view. Horner's rule is
one of them. It is usually presented as an efficient way of evaluating a polynomial, a small
but exact piece of numerical wisdom: fewer multiplications, a better ordering of operations,
a computation that respects the nested structure of the expression. 

In its modern name it is associated with William George Horner
\cite{Horner1819}, although closely related forms had appeared earlier in the
work of Paolo Ruffini \cite{Ruffini1809}, and the underlying idea belongs to a
much older computational tradition \cite{Libbrecht1973}.

Take a polynomial
\[
    p_m(x)=a_0+a_1x+\cdots+a_mx^m.
\]
Horner's form rewrites it as
\[
    p_m(x)=a_0+x\bigl(a_1+x(a_2+\cdots+x(a_{m-1}+xa_m)\cdots)\bigr).
\]
For a given input $x$, equivalently define 
\begin{equation}
    z_0=a_m, \quad z_k= x z_{k-1}+a_{m-k},\quad k=1,\ldots,m.
    \label{eq:horner}
\end{equation}
Then the output is
\[
    p_m(x)=z_m.
\]
This is still just polynomial evaluation. The recurrence changes, however, what one sees. An input is passed through a chain of elementary operations; parameters enter stage by stage; intermediate quantities are propagated forward; the value of the function appears at the end.
\begin{figure}[t]
\centering
\begin{tikzpicture}[>=Latex,node distance=0.72cm]
\tikzstyle{state}=[circle,draw,minimum size=8mm,fill=blue!10]
\tikzstyle{op}=[rectangle,draw,rounded corners,minimum height=7mm,minimum width=22mm,fill=orange!10]
\tikzstyle{ini}=[circle,draw,minimum size=8mm,fill=green!10]

\node[state] (z0) {$z_0$};
\node[op,below=of z0] (op1) {$x z_0 +a_{m-1}$};
\node[state,below=of op1] (z1) {$z_1$};
\node[below=0.55cm of z1] (dots) {$\vdots$};
\node[op,below=0.55cm of dots] (opm) {$x z_{m-1}+a_0$};
\node[state,right=1cm of opm] (zm) {$z_m$};

\node[ini,left=1.5cm of z1] (x) {$x$};

\draw[->,thick] (z0) -- (op1);
\draw[->,thick] (op1) -- (z1);
\draw[->,thick] (z1) -- (dots);
\draw[->,thick] (dots) -- (opm);
\draw[->,thick] (opm) -- (zm);
\draw[->,thick,dashed] (x) to [bend left=20] (op1.west);
\draw[->,thick,dashed] (x) to [bend right=20] (opm.west);

\node[above=0.2cm of zm] {$p_m(x)=z_m$};

\end{tikzpicture}
\caption{\it Classical Horner's rule written as a forward computation.
The input $x$ is shared
across the updates, while the coefficients enter as local parameters
in the recursive rule $z_k=x z_{k-1}+a_{m-k}$.}
\label{fig:horner}
\end{figure}
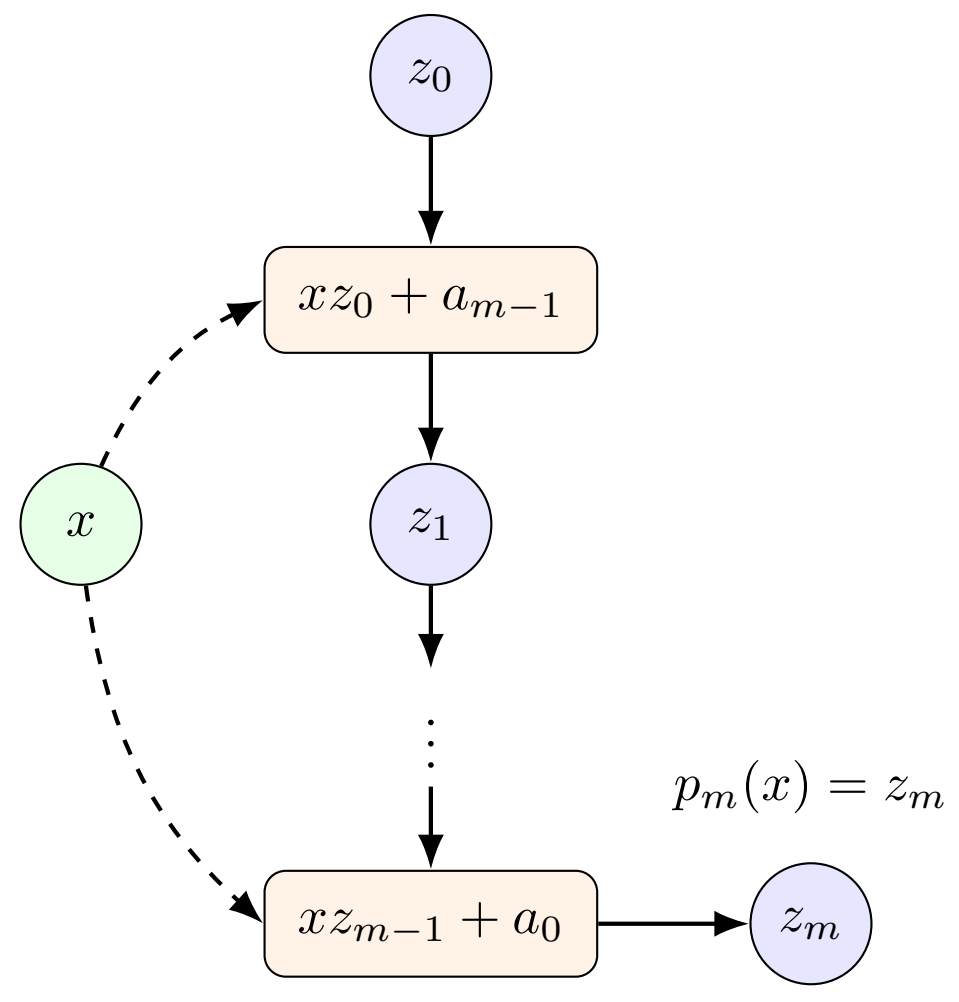

This distinction is worth pausing on. A polynomial can be read as an algebraic expression,
but it can also be read as a computational graph. In one variable this graph is a chain, and the chain
already has a forward pass. Intermediate quantities are produced stage by stage, and the final
value appears at the end of the computation (see Figure \ref{fig:horner}).

The computational-graph viewpoint becomes even more transparent in several variables, where the evaluation is no longer just a chain.
Consider, for example, a bivariate polynomial
\[
p_{m,n}(x_1,x_2)=\sum_{i=0}^{m}\sum_{j=0}^{n} a_{ij}x_1^i x_2^j.
\]
Grouping first in the variable $x_1$, we may write
\[
\begin{split}
p_{m,n}(x_1,x_2)&=q_0(x_2)\\
+ & x_1\bigl(q_1(x_2)+x_1(q_2(x_2)+\cdots+x_1q_m(x_2))\cdots\bigr),
\end{split}
\]
where
\[
q_i(x_2)=a_{i0}+a_{i1}x_2+\cdots+a_{in}x_2^n,\qquad i=0,\dots,m.
\]
Each $q_i(x_2)$ may itself be evaluated by Horner's rule in the variable $x_2$ using 
a chain with intermediate states 
\[
w_{i,j}=w_{i,j-1} x_2 + a_{i,n-j}, \quad j=1,\ldots,n,
\]
starting from $w_{i,0}=a_{i,n}$ to get $
q_i(x_2) = w_{i,n}$.

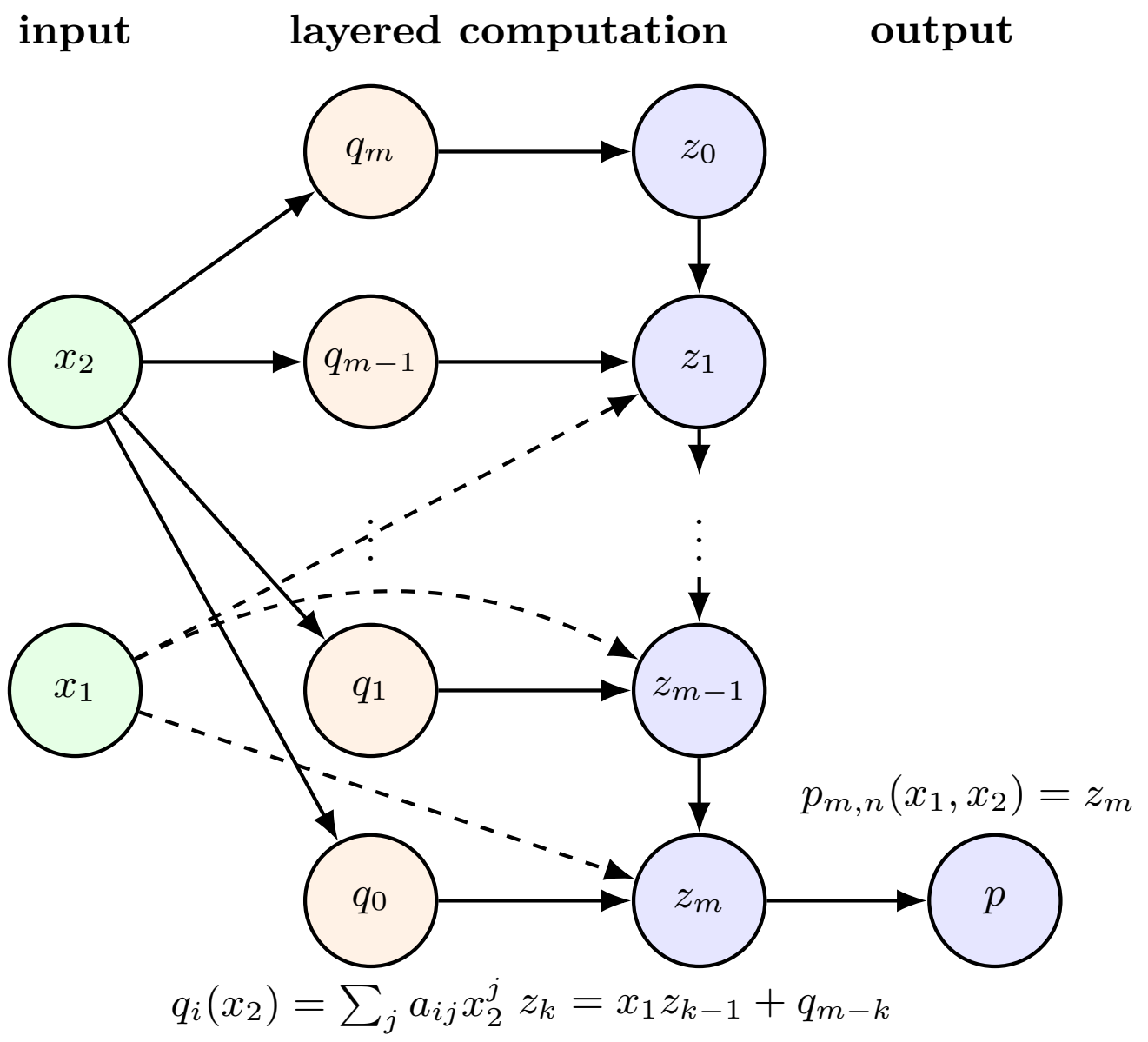
\begin{figure}[tb]
\centering
\begin{tikzpicture}[
    x=1cm,y=1cm,
    >=Latex,
    neuron/.style={circle,draw,thick,minimum size=10mm,inner sep=0pt,font=\small,fill=blue!10},
    neuron2/.style={circle,draw,thick,minimum size=10mm,inner sep=0pt,font=\small,fill=orange!10},
    ini/.style={circle,draw,thick,minimum size=10mm,inner sep=0pt,font=\small,fill=green!10},
    layerlabel/.style={font=\small\bfseries},
    explabel/.style={font=\small},
    conn/.style={->,thick},
    xconn/.style={->,thick,dashed}
]

\node[ini] (yin) at (0,1.5) {$x_2$};
\node[ini] (xin) at (0,-1) {$x_1$};

\node[neuron2] (qm)  at (2.25, 3.1) {$q_m$};
\node[neuron2] (qm1) at (2.25, 1.5) {$q_{m-1}$};
\node[explabel] (qvdots2) at (2.25, 0.25) {$\vdots$};
\node[neuron2] (q1)  at (2.25,-1.0) {$q_{1}$};
\node[neuron2] (q0)  at (2.25,-2.6) {$q_{0}$};

\node[neuron] (hm)  at (4.75, 3.1) {$z_0$};
\node[neuron] (hm1) at (4.75, 1.5) {$z_{1}$};
\node[explabel] (hvdots) at (4.75, 0.25) {$\vdots$};
\node[neuron] (h1)  at (4.75,-1.0) {$z_{m-1}$};
\node[neuron] (h0)  at (4.75,-2.6) {$z_m$};

\node[neuron] (out) at (7,-2.6) {$p$};

\draw[conn] (yin) -- (qm);
\draw[conn] (yin) -- (qm1);
\draw[conn] (yin) -- (q1);
\draw[conn] (yin) -- (q0);

\draw[conn] (qm)  -- (hm);
\draw[conn] (qm1) -- (hm1);
\draw[conn] (q1)  -- (h1);
\draw[conn] (q0)  -- (h0);

\draw[conn] (hm)  -- (hm1);
\draw[conn] (hm1) -- (hvdots);
\draw[conn] (hvdots) -- (h1);
\draw[conn] (h1)  -- (h0);

\draw[xconn] (xin) to[bend left=0] (hm1);
\draw[xconn] (xin) to[bend left=27] (h1);
\draw[xconn] (xin) to[bend left=0] (h0);

\draw[conn] (h0) -- (out);

\node[layerlabel] at (0,4) {input};
\node[layerlabel] at (3.3,4) {layered computation};
\node[layerlabel] at (6.6,4) {output};

\node[explabel,align=center] at (2,-3.4)
{$q_i(x_2)=\sum_{j} a_{ij}x_2^j$};
\node[explabel,align=center] at (4.8,-3.4)
{$z_k=x_1 z_{k-1}+q_{m-k}$};
\node[explabel,align=center] at (6.8,-1.8)
{$p_{m,n}(x_1,x_2)=z_m$};

\end{tikzpicture}
\caption{\it 
A bivariate Horner scheme written as a layered feed-forward computational architecture.
Each node $q_i$ in the first level is itself the output of a one-dimensional Horner
chain in the variable $x_2$. These intermediate quantities are then combined by a
second Horner recursion in the variable $x_1$,
$z_k=x_1 z_{k-1}+q_{m-k}$, until the final output
$z_m=p_{m,n}(x_1,x_2)$ is obtained.}
\label{fig:horner-network}
\end{figure}

If we then define
\[
z_0=q_m(x_2),\quad z_k=x_1\,z_{k-1}+q_{m-k}(x_2),\quad k=1,\dots,m,
\]
we obtain
\[
p_{m,n}(x_1,x_2)=z_m.
\]
Written in this way, the evaluation acquires a layered structure. The variable $x_2$
feeds a first level of computations: each node $q_i(x_2)$ is itself the output of a
one-dimensional Horner chain in $x_2$. These quantities are then passed to a second
stage, where they are combined through an outer Horner recursion in the variable $x_1$.
The variable $x_1$ therefore acts as a shared input to this second stage, and the final
output is the scalar value $p_{m,n}(x_1,x_2)=z_m$.
The resulting architecture is sparse and highly structured, but both
visually and conceptually it already resembles a feed-forward network
(see Figure \ref{fig:horner-network}).

For multivariate polynomials, this architecture is clearly not unique. Different choices of variable
ordering and different ways of factoring and expanding subexpressions lead to different computational graphs,
with different arithmetic costs and, potentially, different numerical properties
\cite{CzekanskySauer2015,KuipersVermaserenPlaatVandenHerik2013}.

I am not identifying a multivariate Horner scheme with a neural network. Its nodes compute algebraic subexpressions rather than learned affine maps followed by nonlinear activations, and its architecture comes from algebraic factorisation rather than from training. The useful observation is more modest: the same polynomial can be represented as a layered computation with inputs, parameters, nodes, intermediate quantities, forward propagation and an output. This gives some of the terminology of machine learning a concrete mathematical starting point.

\subsection{Fitting in a fixed representation}
The next step is to pass from evaluating a function to fitting one. We are given data
\[
    \data=\{(x_i,y_i)\}_{i=1}^n
\]
and choose functions
\[
    \phi_0,\phi_1,\ldots,\phi_m.
\]
For instance, in the simplest polynomial case, $\phi_j(x)=x^j$. We approximate the data
by
\[
    p_m(x)=\sum_{j=0}^m a_j\phi_j(x).
\]
The unknowns are the coefficients
\[
    a=(a_0,\ldots,a_m)^T.
\]
Evaluating the basis functions at the data points gives the matrix
\[
    V_{ij}=\phi_j(x_i),\qquad i=1,\ldots,n,\quad j=0,\ldots,m.
\]
Thus the fitting problem is
\begin{equation}
    Va\simeq y,
    \label{eq:fit}
\end{equation}
where $y=(y_1,\ldots,y_n)^T$.

This is the canonical fixed-representation problem, and it is one of the basic
viewpoints of classical approximation theory and numerical linear algebra
\cite{Trefethen2013,GolubVanLoan2013}. The functions $\phi_j$ are chosen in
advance. The data determine the rows of $V$. The coefficients are then selected by an
error criterion. In the Euclidean norm we obtain
\begin{equation}
    \min_{a}\|Va-y\|_2^2.
    \label{eq:mina}
\end{equation}

At this level, supervised learning has not yet introduced a new mathematical object. The
training set is the data set. The model class is the approximation space. The loss is a
measure of the residual. The training algorithm is the numerical method used to determine
the parameters. This translation is simple, but it is not superficial. It is the point at
which the language of numerical analysis begins to organise the language of learning.

\subsection{Three regimes in one matrix equation}

Let
\[
    p=m+1
\]
be the number of coefficients. Equation \eqref{eq:fit} contains three regimes that are often taught
separately:
\[
    p<n,\qquad p=n,\qquad p>n.
\]
The classification above is algebraic: it concerns the size and rank structure of the
matrix problem. Its interpretation also depends on the nature of the data. If the values
$y_i$ are exact samples of an underlying function, the square case naturally leads to
interpolation, and the overparametrised case to many exact interpolants (see Table \ref{tab:fitting-regimes}). If, instead, the
data contain measurement error, sampling error, or unresolved modelling effects, exact
fitting has a different meaning: it may reproduce the error rather than the underlying
quantity of interest. In that case the residual is not merely a computational defect, but
part of the modelling problem. This is why losses, regularisation and validation enter
naturally in the next step.

\begin{table}[htb]
\centering
\small
\setlength{\tabcolsep}{3pt}
\begin{tabular}{ccl}
\toprule
Regime & Matrix problem & Numerical meaning \\
\midrule
$p<n$ & overdetermined & least-squares approximation \\
$p=n$ & square & interpolation ($V$ nonsingular) \\
$p>n$ & underdetermined & many possible interpolants \\
\bottomrule
\end{tabular}
\caption{\emph{The different regimes of a fitting problem.}}
\label{tab:fitting-regimes}
\end{table}

In the first case, exact fitting is generally impossible, and the residual must be
balanced. If $V$ has full column rank, the least-squares solution is unique and satisfies
\[
    V^TVa=V^Ty.
\]
A first course may introduce the normal equations. A numerical analyst immediately asks how the computation should actually be performed, and replaces them when necessary by QR or SVD factorisations \cite{GolubVanLoan2013}. A useful didactic point is often hidden here. In elementary presentations, the overdetermined case is usually introduced through the regression line: many data points, two unknown coefficients, and no exact line passing through all observations. This example is important, but it can narrow the discussion too much. What matters is not the straight line itself; it is approximation in a prescribed finite-dimensional space.

The same matrix equation covers the fitting of polynomials, trigonometric expansions, splines, radial basis functions, finite elements, or any chosen family of trial functions. What changes is the representation, not the basic structure of the problem. The regression line is therefore not a separate elementary topic. It is the first visible instance of a much more general approximation problem.

In the second case, for instance with the monomial basis and distinct nodes, $V$ is a
Vandermonde matrix and interpolation becomes possible. This is the familiar theorem of
existence and uniqueness. But numerical analysis adds the warning that uniqueness is not
stability. Interpolation may fit the data exactly and still behave badly between the data
points or under perturbations.

In the third case, the system may have infinitely many exact solutions. If $Va=y$ is
consistent and $V$ has full row rank, the minimum Euclidean norm interpolant is
\begin{equation}
    a^\dagger=V^\dagger y=V^T(VV^T)^{-1}y.
    \label{eq:3}
\end{equation}
This formula is elementary, but its interpretation matters. Overparametrisation does not first appear in deep learning; it appears already in linear algebra as the loss of uniqueness. The mathematical question is then not whether the data can be fitted, but which fitting solution is selected \cite{BelkinHsuXu2020,HastieMontanariRossetTibshirani2022}.

Selection may come from a norm, as in \eqref{eq:3}. It may come from regularisation, for example
\[
    \min_a \|Va-y\|_2^2+\lambda\|a\|_2^2,
\]
which gives, for $\lambda>0$,
\begin{equation}
    a_\lambda=(V^TV+\lambda I)^{-1}V^Ty.
\end{equation}
It may come from an iterative method, from a stopping rule, or from the parametrisation itself. Here numerical analysis and modern machine learning begin to share a language. The perspective followed here is deliberately numerical. It does not aim to give a genealogy of machine learning, nor to cover its statistical and probabilistic foundations. Instead, it uses the familiar fitting problem $Va\simeq y$ as an entry point: the residual becomes a loss, the solution process becomes training, and the numerical method used to reduce the loss becomes the optimiser.

\section{Losses, gradients and optimisation}

\subsection{The error is part of the model}

The Euclidean norm is the natural entry point, but it should not be mistaken for a neutral
choice. The problem \eqref{eq:mina}
penalises large residuals strongly and leads to projection in Euclidean geometry. The
problem
\[
    \min_a\|Va-y\|_1
\]
responds differently to outliers. The problem
\[
    \min_a\|Va-y\|_\infty
\]
asks for control of the worst residual and leads toward minimax approximation. Thus the loss is not an implementation detail. It is a modelling choice.
Changing the loss changes the geometry of the problem and the numerical tools
available for solving it \cite{Bishop2006}.

When the loss is quadratic and the parametrisation is linear, direct methods are often the
natural numerical tools. The gradient viewpoint, however, remains useful even in this
classical case, because it exhibits the loss as an object to be minimised and the rows of
$V$ as the elementary contributions to that minimisation. When the loss is
non-quadratic, non-smooth, large-scale, or nonlinear in the parameters, iterative
optimisation is no longer merely an alternative implementation. It becomes part of the
mathematical description of the computed model.

\subsection{The loss and the computation of its gradient}
The same fitting problem also explains why optimisation enters naturally. Even before
neural networks are introduced, the coefficients may be selected iteratively by minimising
a loss. Write
\begin{equation}
    \mathcal L(a)
     = 
    \frac{1}{2n}
    \sum_{i=1}^n
    \bigl(p_m(x_i;a)-y_i\bigr)^2 = \frac1{2n}\|Va-y\|_2^2.
\end{equation}
If
$
    r_i(a)=p_m(x_i;a)-y_i
$
is the residual at the data point $x_i$, then differentiating the loss with respect to
the coefficient $a_j$ gives
\[
    \frac{\partial \mathcal L}{\partial a_j}(a)
    =
    \frac{1}{n}
    \sum_{i=1}^n
    r_i(a)
    \frac{\partial p_m}{\partial a_j}(x_i;a).
\]
For a fixed representation
\[
    p_m(x;a)=\sum_{j=0}^m a_j\phi_j(x),
\]
one has
\[
    \frac{\partial p_m}{\partial a_j}(x_i;a)=\phi_j(x_i).
\]
Therefore
\[
    \frac{\partial \mathcal L}{\partial a_j}(a)
    =
    \frac{1}{n}
    \sum_{i=1}^n
    r_i(a)\phi_j(x_i),
    \qquad j=0,\ldots,m.
\]
Equivalently,
\[
    \nabla_a \mathcal L(a)
    =
    \frac{1}{n}V^T(Va-y).
\]
A gradient descent iteration therefore has the form
\begin{equation}
    a^{r+1}
    =
    a^r-\eta\,\nabla_a\mathcal L(a^r)
    =
    a^r-\eta\,\frac{1}{n}V^T(Va^r-y),
\end{equation}
where $\eta$ is the step size.

If $B\subset\{1,\ldots,n\}$ is a subset of the data points, chosen at random, the corresponding
mini-batch loss is
\begin{equation}
    \mathcal L_B(a)
    =
    \frac{1}{2|B|}
    \sum_{i\in B}
    \bigl(p_m(x_i;a)-y_i\bigr)^2,
\end{equation}
and the mini-batch gradient is
\[
    \frac{\partial \mathcal L_B}{\partial a_j}(a)
    =
    \frac{1}{|B|}
    \sum_{i\in B}
    r_i(a)\phi_j(x_i),
\]
which reduces to $ |B|^{-1}\sum_{i\in B} r_i(a)x_i^j $ in the monomial basis.
Thus
\[
    a^{r+1}
    =
    a^r-\eta\,\nabla_a\mathcal L_B(a^r)
\]
is already the mini-batch stochastic gradient method in the
elementary polynomial fitting problem \cite{BottouCurtisNocedal2018}.

In the monomial case $\phi_j(x)=x^j$, the same calculation can be written without expanding the powers $x_i^j$. For one data
point $(x,y)$, evaluate the polynomial by Horner's recurrence
\[
    z_0=a_m,\qquad
    z_k=xz_{k-1}+a_{m-k},\quad k=1,\ldots,m,
\]
so that
\[
    p_m(x;a)=z_m.
\]
For the single-data loss
\[
    \mathcal L(a)=\frac12\bigl(z_m-y\bigr)^2,
\]
set
\[
    \bar z_m
    =
    \nabla_{z_m}\mathcal L
    =
    z_m-y.
\]
The chain rule applied backward to the recurrence gives, for
$k=m,m-1,\ldots,1$,
\[
    \frac{\partial \mathcal L}{\partial a_{m-k}}
    =
    \bar z_k,
    \qquad
    \bar z_{k-1}
    =
    x\bar z_k.
\]
At the end of the backward recursion, since $z_0=a_m$, one also has
\[
    \frac{\partial \mathcal L}{\partial a_m}
    =
    \bar z_0.
\]
For a mini-batch $B$, the same idea can be written directly on the batch loss starting from the recurrence 
\[
 z_{i,0}=a_m,\qquad
    z_{i,k}=x_i z_{i,k-1}+a_{m-k},\quad k=1,\ldots,m,
\]
to get 
\[
\frac{\partial \mathcal L_B}{\partial a_m}
=\frac{1}{|B|} \sum_{i\in B}\bar z_{i,0}.
\]
Thus the gradients with respect to all coefficients are obtained by one
forward evaluation on the batch followed by one backward recursion through
the same Horner graph for each data point, with the resulting contributions
averaged over the batch. Explicitly,
\[
    \nabla_a \mathcal L_B(a)
    =
    \frac{1}{|B|}\sum_{i\in B}
    \bigl(z_{i,m}-y_i\bigr)
    \bigl(1,x_i,\ldots,x_i^m\bigr)^T .
\]
This is a simple prototype of backpropagation. In the polynomial case, the
computational graph is just Horner's chain and the parameters are the coefficients
$a_j$. In a neural network, the graph is larger and the parameters are the weights and
biases, but the principle is the same: a loss is differentiated with respect to many
parameters by propagating adjoint variables backward through a forward computation.

\begin{figure}[tb]
\begin{center}
\includegraphics[width=0.49\textwidth]{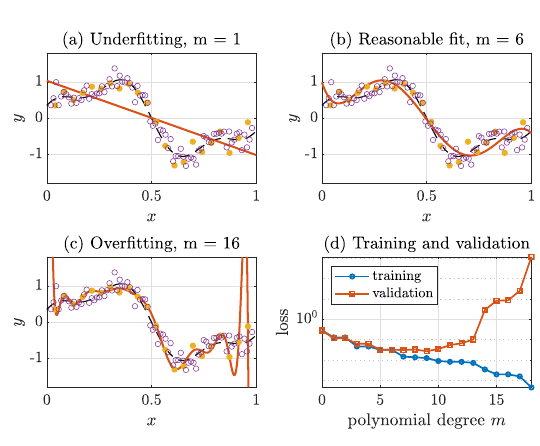}
\end{center}
\vskip -.35cm
\caption{\it
Polynomial fitting of noisy data exhibits the basic phenomena of learning. The orange points and the purple circles
denote the training and validation data.
Panels (a)--(c) show polynomial fits of increasing degree:
underfitting, a reasonable fit, and overfitting. Panel (d) shows the corresponding
training and validation losses as functions of the degree. As the degree increases, the training loss decreases, while the validation loss may
increase once the polynomial starts fitting fluctuations in the training data.
}
\label{fig:poly-overfitting}
\end{figure}

A further distinction is essential in learning. The optimiser sees only the training loss: the loss computed on the data used to select the parameters. After the parameters have been selected, the same approximation can be evaluated on data that were not used in the optimisation. If we denote by $\{\tilde x_i, \tilde y_i\}_{i=1}^r$ 
a validation set and by $\tilde V$ the corresponding matrix, one may compare
\[
    {\mathcal L}(a)
    = \frac1{2n} \|Va-y\|^2_2,
\]
with
\[
    \tilde {\mathcal L}(a)
    = \frac1{2r} \|\tilde Va-\tilde y\|^2_2.
\]
The first quantity is the objective seen by the optimiser. The second tests how well the selected approximation generalises beyond the data that selected it. This distinction is already meaningful for polynomial fitting. Increasing the
degree enriches the approximation space and can reduce the training loss. But if the data
are affected by noise, the same process may eventually increase the validation loss: the
polynomial begins to fit the fluctuations of the training data rather than the underlying
function. This is the classical manifestation of overfitting, and it is shown in
Figure~\ref{fig:poly-overfitting}.

At this point the vocabulary can be translated without ceremony. Words such as training, features, loss, optimiser, batch or epoch may suggest that a new mathematical universe has appeared. In the polynomial problem, these words already have transparent meanings (see Table \ref{tab:ml-polynomial-dictionary}).

\begin{table*}[t]
\centering
\small
\renewcommand{\arraystretch}{1.16}
\begin{tabularx}{\textwidth}{p{0.21\textwidth}p{0.30\textwidth}Y}
\toprule
\textbf{Machine learning} &
\textbf{Polynomial fitting} &
\textbf{Numerical meaning} \\
\midrule
training set &
$\{(x_i,y_i)\}_{i=1}^n$ &
data set used to determine the coefficients $a_j$ \\

validation set &
$\{(\tilde x_i,\tilde y_i)\}_{i=1}^r$  &
data not used to determine the coefficients $a_j$\\

model &
$p_m(x)=\sum_{j=0}^m a_j\phi_j(x)$ &
approximation space \\

features &
$\phi_j(x_i)$ &
basis functions evaluated on data \\

parameters &
$a_0,\ldots,a_m$ &
unknown coefficients \\

prediction &
$p_m(x_i)$ &
computed approximation \\

loss &
$\mathcal L(a)=\frac1{2n}\|Va-y\|_2^2$ &
residual measure \\

training &
$\min_a \|Va-y\|_2^2$ &
parameter identification \\

optimiser &
QR, SVD, gradient descent, SGD &
numerical method used to select parameters \\

learning rate &
$\eta$ in $a^{r+1}=a^r-\eta\nabla\loss(a^r)$ &
step size of an iterative method \\

epoch &
one pass through the data &
one complete pass over the rows of $V$ \\

batch or mini-batch &
$B\subset\{1,\ldots,n\}$, with rows $V_B$ &
subset of data used to form a residual or gradient \\
regularisation &
$\lambda\|a\|_2^2$ &
stabilised selection \\

generalisation &
$\tilde {\mathcal L}(a)= \frac1{2r} \|\tilde Va -\tilde y\|_2^2$ &
accuracy away from data \\
\bottomrule
\end{tabularx}
\caption{\it A small dictionary between the vocabulary of machine learning and the
classical language of polynomial fitting.}
\label{tab:ml-polynomial-dictionary}
\end{table*}

\section{Learning the approximation space}

\subsection{From prescribed to parametrised representations}
The polynomial problem fixes the representation. The functions $\phi_j$ are prescribed,
and the matrix $V$ follows from them. The data determine the rows, but the columns are
chosen before the fitting begins.

Many learning models change this arrangement: the representation itself is parametrised. Instead of fixed functions $\phi_j$, one has
\[
    \phi_j=\phi_j(\cdot;\theta),
\]
and the empirical feature matrix becomes
\[
    V(\theta)_{ij}=\phi_j(x_i;\theta).
\]
The fitting problem then has the form
\begin{equation}
    \min_{\theta,a}\|V(\theta)a-y\|_2^2.
    \label{eq:fit2}
\end{equation}
Equation \eqref{eq:fit2} is not intended as standard machine-learning notation. Its value here is to separate two operations that are often combined in a single parameter vector. Some parameters determine the representation; others select an element in the representation so obtained. In many neural networks these roles are trained simultaneously, but the distinction is useful pedagogically and numerically.
For fixed $\theta$, one is back to a fitting problem in a prescribed empirical space. Varying $\theta$ changes that space \cite{Bishop2006,GoodfellowBengioCourville2016,HighamHigham2019}.

The need to change a representation does not first arise with neural networks. It is already visible in classical approximation. A fixed basis may be poorly adapted to the data, badly conditioned, or too large to be useful in high dimension. Polynomial approximation gives a familiar example. High-degree global polynomials are expressive, but they may oscillate, become ill-conditioned, and propagate local perturbations through a global representation.

One classical response is to change the structure of the approximation space. Splines, for instance, replace a single global polynomial by a piecewise polynomial representation,
where the functions $\phi_j$ are local basis functions determined by the degree and by the chosen knots. The fitting problem again 
has the form \eqref{eq:fit} but the representation is now local and prescribed rather than global and prescribed. This locality is one reason for the good numerical behaviour of splines: increasing the number of knots enriches the space without increasing one polynomial degree everywhere.

This distinction is useful when moving toward learning. Polynomial bases and spline bases are different representations, but in both cases the representation is chosen before the fitting begins. The data select the coefficients, not the basis functions themselves. In neural networks, by contrast, the representation is trainable.

This suggests a broader approximation problem. Instead of choosing only the coefficients in
\begin{equation}
    f_a(x)=\sum_{j=1}^p a_j\phi_j(x),
    \label{eq:class}
\end{equation}
one also allows the functions themselves to depend on parameters:
\begin{equation}
    f_{\theta,a}(x)
    =
    \sum_{j=1}^p a_j\phi_j(x;\theta).
\end{equation}
The parameters $\theta$ may describe centres and scales in radial basis approximation,
knots in an adaptive spline, frequencies in an exponential expansion, atoms in a dictionary, or latent factors in a low-rank model. In all these
cases, for fixed $\theta$ the problem is again a linear fitting problem for the
coefficients $a$. Varying $\theta$, however, changes the columns of the empirical
matrix and hence changes the approximation space itself.

This is the bridge from approximation to learning used here. Some parameters determine the representation in which the final fitting problem is posed. Neural networks are one important realisation of this idea, but not its first appearance.

\subsection{Neural networks as adaptive approximation spaces}

Let us return for a moment to the classical approximation problem \eqref{eq:class}. 
It is useful to collect the functions $\phi_j$ into the vector
\[
    \Phi(x)=
    \bigl(\phi_1(x),\ldots,\phi_p(x)\bigr)^T.
\]
The approximant can then be written as
\[
    f_a(x)=a^T\Phi(x).
\]
In this notation, a feature is simply one of the quantities $\phi_j(x)$ used in
the final linear combination. The map $\Phi$ collects these quantities. In
polynomial approximation, for example, one may take
\[
    \Phi_m(x)=
    \bigl(1,x,x^2,\ldots,x^m\bigr)^T,
\]
or, in several variables, a vector containing all monomials up to a prescribed
degree. Horner's rule showed another viewpoint on the same polynomial object: the final value can be obtained by propagating intermediate states through a computation, rather than by forming all powers separately. 

A feed-forward neural network uses this computational viewpoint at the level of the representation itself. It keeps the idea of a final linear combination, but constructs the vector that enters this combination by a sequence of state updates.
Let $x\in\mathbb{R}^d$ be the input and set
\[
    z_0(x)=x.
\]
The successive layers define states
\begin{equation}
z_k(x)=\sigma\bigl(W_k z_{k-1}(x)+b_k\bigr),
\qquad k=1,\ldots,L.
\end{equation}
Here
\[
W_k\in\mathbb{R}^{n_k\times n_{k-1}},
\qquad
b_k\in\mathbb{R}^{n_k},
\qquad n_0=d.
\]
The matrix $W_k$ contains the weights of the $k$-th layer: it determines how
the quantities produced by the previous layer are linearly combined. The vector
$b_k$ contains the biases: it shifts these combinations before the scalar
function is applied. 

The function $\sigma$, called the activation function, is
applied componentwise. 
Its role is part of the numerical structure of the network: it affects
both the class of representations that can be produced and the stability of
their propagation through the layers. 

The last state $z_L(x)$ is the vector that will be linearly combined. We write
\[
    \Phi_\theta(x)=z_L(x)
    =
    \bigl(\phi_1(x;\theta),\ldots,\phi_p(x;\theta)\bigr)^T,
    \qquad p=n_L,
\]
where
\[
    \theta=\{W_k,b_k\}_{k=1}^L
\]
denotes the parameters used to construct the layers. The approximating function is
then
\[
    F_{\theta,a}(x)
    =
    a^T\Phi_\theta(x)
    =
    \sum_{j=1}^p a_j\phi_j(x;\theta).
\]
Thus the final form is the same as in classical approximation: a linear combination
of functions. What has changed is the origin of those functions \cite{Pinkus1999,Barron1993}.

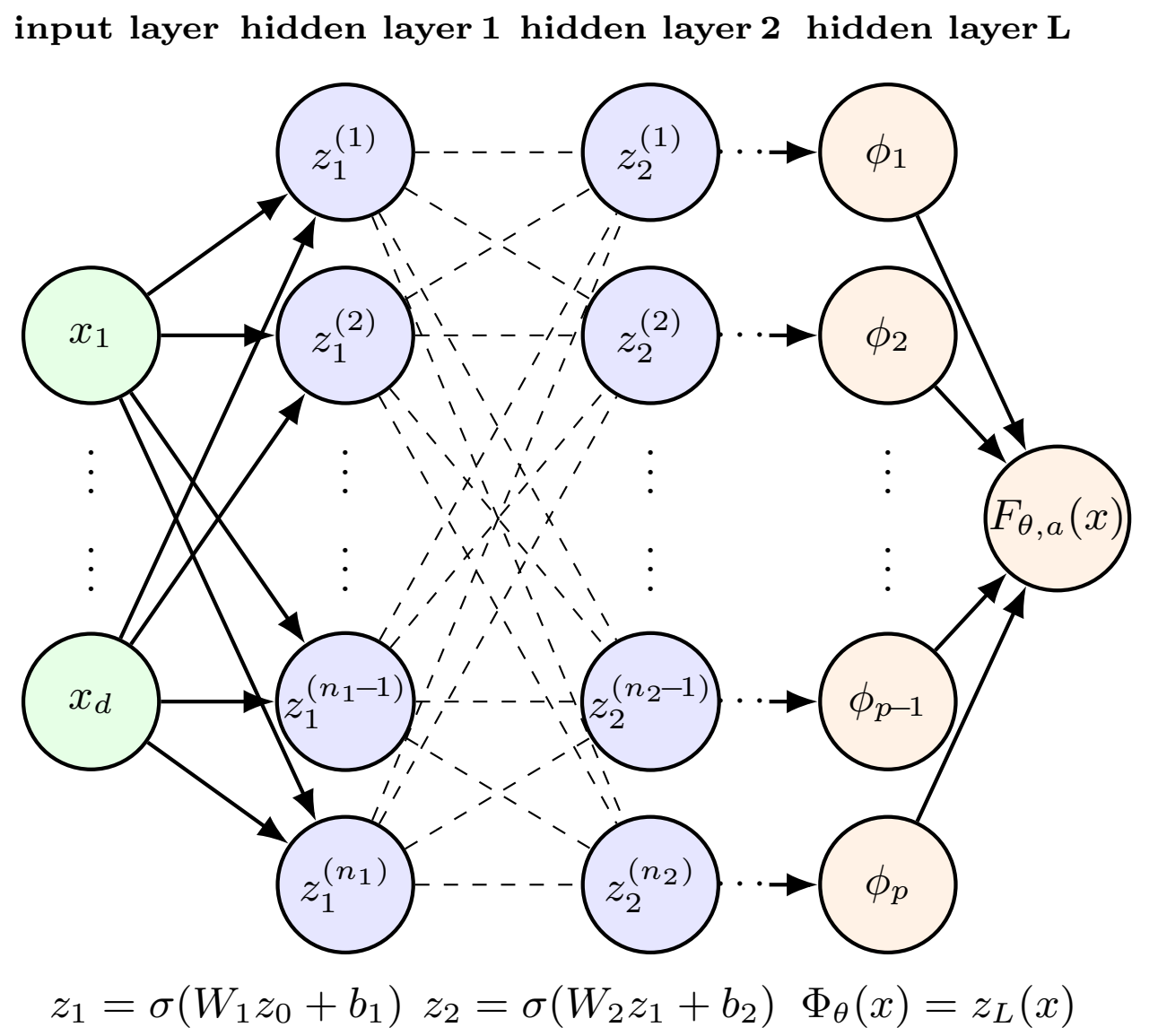
\begin{figure}[tb]
\centering
\begin{tikzpicture}[
    >=Latex,
    x=1.25cm,
    y=0.9cm,
    state/.style={circle,draw,thick,minimum size=10mm,inner sep=0pt,font=\small,fill=blue!10},
    ini/.style={circle,draw,thick,minimum size=10mm,inner sep=0pt,font=\small,fill=green!10},
    state2/.style={circle,draw,thick,minimum size=10mm,inner sep=0pt,font=\small,fill=orange!10},
    out/.style={rectangle,draw,rounded corners,minimum height=7mm,minimum width=17mm},
    every node/.style={font=\small}
]

\node[ini] (xone) at (-0.5,0.5) {$x_1$};
\node at (-0.5,-.5) {$\vdots$};
\node at (-0.5,-1.3) {$\vdots$};
\node[ini] (xd) at (-0.5,-2.5) {$x_d$};

\node[state] (z11) at (1,2) {$z_1^{(1)}$};
\node[state] (z12) at (1,0.5) {$z_1^{(2)}$};
\node at (1,-0.5) {$\vdots$};
\node at (1,-1.3) {$\vdots$};
\node[state] (z1n1) at (1,-2.5) {$z_1^{(n_1\!-\!1)}$};
\node[state] (z1n) at (1,-4) {$z_1^{(n_1)}$};

\node[state] (z21) at (2.8,2) {$z_2^{(1)}$};
\node[state] (z22) at (2.8,0.5) {$z_2^{(2)}$};
\node at (2.8,-0.5) {$\vdots$};
\node at (2.8,-1.3) {$\vdots$};
\node[state] (z2n1) at (2.8,-2.5) {$z_2^{(n_2\!-\!1)}$};
\node[state] (z2n) at (2.8,-4) {$z_2^{(n_2)}$};

\node at (3.35,2) {$\cdots$};
\node at (3.35,0.5) {$\cdots$};
\node at (3.35,-2.5) {$\cdots$};
\node at (3.35,-4) {$\cdots$};

\node[state2] (zL1) at (4.2,2) {$\phi_1$};
\node[state2] (zL2) at (4.2,0.5) {$\phi_2$};
\node at (4.2,-0.5) {$\vdots$};
\node at (4.2,-1.3) {$\vdots$};
\node[state2] (zLp1) at (4.2,-2.5) {$\phi_{p\!-\!1}$};
\node[state2] (zLp) at (4.2,-4) {$\phi_p$};
\node at (4.5,-5) {$\Phi_\theta(x)=z_L(x)$};

\node[state2] (out) at (5.2,-1) {$F_{\theta,a}(x)$};

\foreach \u in {xone,xd}
  \foreach \v in {z11,z12,z1n1,z1n}
    \draw[->,thick] (\u) -- (\v);

\foreach \u in {z11,z12,z1n1,z1n}
  \foreach \v in {z21,z22,z2n1,z2n}
    \draw[dashed] (\u) -- (\v);

\draw[->,thick] (3.5,2) -- (zL1);
\draw[->,thick] (3.5,0.5) -- (zL2);
\draw[->,thick] (3.5,-2.5) -- (zLp1);
\draw[->,thick] (3.5,-4) -- (zLp);

\draw[->,thick] (zL1) -- (out);
\draw[->,thick] (zL2) -- (out);
\draw[->,thick] (zLp1) -- (out);
\draw[->,thick] (zLp) -- (out);

\node at (-.35,3.0) {\textbf{\scriptsize input layer}};
\node at (.3,-5) {$z_1=\sigma(W_1z_0+b_1)$};
\node at (2.5,-5) {$z_2=\sigma(W_2z_1+b_2)$};

\node at (1.15,3.0) {\textbf{\scriptsize hidden layer\! 1}};
\node at (2.8,3.0) {\textbf{\scriptsize hidden layer\! 2}};
\node at (4.5,3.0) {\textbf{\scriptsize hidden layer\! L}};

\end{tikzpicture}
\caption{\it A feed-forward neural network.
Each column is a layer, and the scalar components of the layer are its nodes.
The hidden layers construct the vector $\Phi_\theta(x)=z_L(x)$. The final
operation is the linear combination $F_{\theta,a}(x)=a^T\Phi_\theta(x)$.}
\label{fig:ffnn-feature-map}
\end{figure}

This gives the usual layered structure of a feed-forward neural network
(see Figure~\ref{fig:ffnn-feature-map}). The depth is the number of successive transformations
applied to the input, while the width of a layer is the number of scalar quantities
carried by that layer. These two notions answer different questions. Depth measures
how many stages of composition are used. Width measures how many intermediate
features are propagated in parallel at a given stage.

A useful way to interpret this distinction is to return for a moment to polynomial approximation. In one variable, increasing the degree adds new powers of $x$ and therefore enriches the approximation space. In several variables, however, a single degree already contains many different features. For example, in two variables the homogeneous terms of degree $k$ are
\[
x_1^k,\quad x_1^{k-1}x_2,\quad \ldots,\quad x_1x_2^{k-1},\quad x_2^k.
\]
Thus the degree may be viewed as a level of algebraic complexity, while the monomials of that degree are parallel features carried at the same level. In $d$ variables, the number of monomials of degree $k$ is
\[
\binom{d+k-1}{k},
\]
which grows rapidly with both $d$ and $k$. This is one elementary form of the curse of dimensionality.

This polynomial picture gives a useful guide to the role of width and depth in neural networks.
A wide layer can carry
many intermediate features in parallel, much as a polynomial space in many variables
contains many monomials of the same degree. Additional layers allow features to be
combined and transformed through successive stages, producing increasingly complex
representations. 
The analogy should not be taken literally: neural networks do not
usually generate all polynomial features of a given degree. Their strength is precisely
that the intermediate features are not fixed in advance, but are adapted by the
parameters $\theta$ and by the optimisation process. Still, the polynomial analogy
explains why both width and depth matter: width controls how many features can be
represented at a given stage, while depth controls how features are successively
composed.

Evaluating the features on the data gives the empirical matrix
\[
    V(\theta)_{ij}
    =
    \phi_j(x_i;\theta)
    =
    \bigl(\Phi_\theta(x_i)\bigr)_j .
\]
Thus, for fixed $\theta$, fitting the final coefficients $a$ is again a linear problem,
\[
    V(\theta)a\simeq y.
\]
But $\theta$ is not fixed during training. Changing $\theta$ changes the functions
$\phi_j(\cdot;\theta)$, hence changes the columns of $V(\theta)$, and therefore changes
the approximation space itself. 

\section{Tasks, training and backpropagation}

\subsection{Changing the task by changing the loss}
The same architecture can be used for different tasks by changing the loss. In the
least-squares setting discussed so far, the output $F_{\theta,a}(x)$ is compared with a
scalar value $y$ through a quadratic residual. 

\begin{figure}[tb]
\begin{center}
\includegraphics[width=0.49\textwidth]{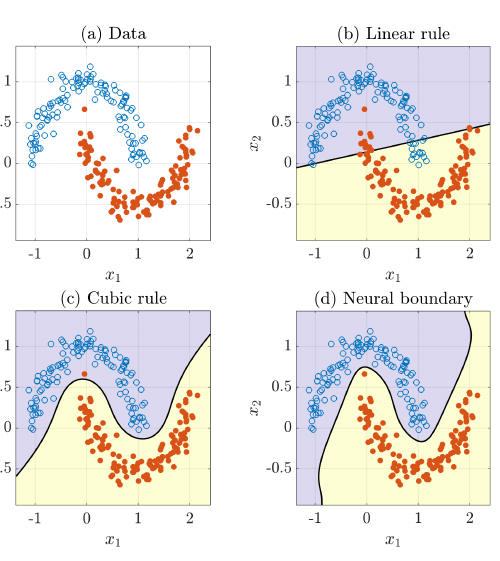}
\end{center}
\vskip -.4cm
\caption{\it
A comparison between prescribed and learned representations in a two-class classification problem. The blue circles and the red points denote the two classes. A linear rule in the original variables gives only a straight decision boundary. A prescribed cubic rule can already produce a nonlinear boundary and improve the separation substantially. In the neural-network case, the representation is not fixed in advance, but is
constructed by training.
}
\label{fig:learned}
\end{figure}

For a two-class classification problem, the data no longer represent values of a real-valued function, but labels $y_i\in\{0,1\}$. One may therefore view the task as the approximation of the indicator of one class. The network first constructs a real-valued function
\[
u_{\theta,a}(x)=a^T\Phi_\theta(x),
\]
where a bias term can be included by adding a constant component to $\Phi_\theta$. Since an indicator takes values between $0$ and $1$, this function is then passed through the logistic map
\[
\pi_{\theta,a}(x)=\frac{1}{1+\exp(-u_{\theta,a}(x))}.
\]
The quantity $\pi_{\theta,a}(x)$ may be interpreted as a smooth approximation of the class indicator. A standard empirical discrepancy for this approximation is the cross-entropy
\[
\mathcal L(\theta,a)
\!=\!
-\frac1n\sum_{i=1}^n
\left[
y_i\log \pi_{\theta,a}(x_i)
\!+\!
(1\!-\!y_i)\log\bigl(1\!-\!\pi_{\theta,a}(x_i)\bigr)
\right].
\]
This loss penalises small values of $\pi_{\theta,a}(x_i)$ on points with $y_i=1$ and large values on points with $y_i=0$. Thus the task has changed from fitting real-valued data to separating two classes, but the representational structure is the same. A fixed representation, for example one built from polynomial features, can also produce nonlinear classifiers in the original variables (see Figure \ref{fig:learned}). The distinction is between prescribing such a representation before fitting and adapting the representation during training.

This makes the architecture task-flexible. Changing the loss changes the task without necessarily changing the basic computational architecture.
In regression the loss measures residuals in observed values; in classification it
measures the mismatch between predicted and observed labels; in physics-informed neural
networks one adds terms measuring the residual of a differential equation and the
mismatch with boundary or initial data \cite{RaissiPerdikarisKarniadakis2019}. Thus data fitting, classification and
physics-constrained approximation can all be expressed as optimisation problems over the
same class of parametrised computational graphs.

\subsection{Differentiating the computational graph}
This also clarifies the connection with Horner's rule. Horner's recurrence
is an ordered computation that propagates intermediate states. A feed-forward neural
network uses the same computational syntax, but replaces the fixed scalar update
\[
    z\mapsto xz+a
\]
by vector-valued layer maps
\[
    z\mapsto \sigma(W_kz+b_k).
\]
The comparison concerns the syntax of the computation--input, states, layerwise propagation and output--not the identity of the elementary maps.

The backpropagation calculation seen for Horner's rule extends directly to this
layered architecture. Once a function is
represented by a sequence of elementary maps,
\[
    z_0\mapsto z_1\mapsto\cdots\mapsto z_L\mapsto F_{\theta,a}(x),
\]
its derivatives are obtained by applying the chain rule along that sequence. 

To display the backward propagation without the additional notation of classification
losses, let us return to the quadratic loss.
Consider, for simplicity, one data point $(x,y)$ and the loss
\[
    \mathcal{L}(\theta,a)
    =
    \frac12\bigl(F_{\theta,a}(x)-y\bigr)^2.
\]
The gradient of the loss with respect to the final
coefficients is
\[
    \nabla_a\mathcal{L} = (F_{\theta,a}(x)-y) z_L.
\]
The adjoint variables, that is the gradients of the loss with respect to the
intermediate states, are then propagated backward by the chain rule. Starting from
\[
    \bar z_L
    =
    \nabla_{z_L}\mathcal{L}
    =
    (F_{\theta,a}(x)-y) a,
\]
the backward recursion for
$k=L,L-1,\ldots,1$, gives
\[
    \delta_k
    =
    \bar z_k \odot \sigma'(W_kz_{k-1}+b_k),
\]
where $\odot$ denotes componentwise multiplication. Then
\[
    \nabla_{W_k}\mathcal{L}
    =
    \delta_k z_{k-1}^T,
    \qquad
    \nabla_{b_k}\mathcal{L}
    =
    \delta_k,
\]
and the adjoint variable is propagated to the previous layer by
\[
    \bar z_{k-1}
    =
    W_k^T\delta_k.
\]
Thus one backward pass produces the gradients with respect to all weights and biases.
For a batch $B$, the corresponding gradients are obtained by averaging these
contributions over $(x_i,y_i)\in B$.

This is the computational advantage. A direct application of the chain rule
separately to each parameter would repeatedly traverse the same computational graph.
Backpropagation instead reuses the intermediate quantities computed in the forward pass
and propagates the gradients with respect to the intermediate states once in reverse order.
The cost of evaluating all gradients is therefore of the same order as a small number of
forward evaluations, rather than the cost of differentiating separately with respect to
each parameter. 
In this sense, backpropagation is not a new principle attached only to neural
networks. It is reverse-mode differentiation of a parametrised
computational graph, organised so that gradient-based optimisation becomes feasible
\cite{GriewankWalther2008,RumelhartHintonWilliams1986}.

The role of the starting point in a gradient descent iteration also changes. 
In a fixed representation, an
initial vector $a^0$ is a guess for the coefficients; the approximation space is
already prescribed. In a trainable representation, the initial values
$(\theta^0,a^0)$ have two roles. The parameter $\theta^0$ determines the
initial feature map $\Phi_{\theta^0}$, and hence the initial empirical matrix
$V(\theta^0)$. The vector $a^0$ determines the initial linear combination of
these features.  For this reason, neural-network training usually starts from random values of
the weights, while biases are often initialised separately, for instance at zero
or at small values. 

The structure of \eqref{eq:fit2} also has a classical numerical analogue. The problem is linear in the coefficients
$a$ and nonlinear in the parameters $\theta$. For fixed $\theta$, the coefficients can
be selected by a linear least-squares problem. Varying $\theta$ changes the space in
which that linear problem is posed. This separable structure is familiar from variable
projection methods for nonlinear least squares \cite{GolubPereyra1973}. Neural-network
training is different in scale, architecture and nonconvexity, but the underlying
separation remains useful.

\section{Putting the pieces together}

\subsection{A numerical perspective on learning}
The route can now be summarised in numerical terms. Horner's rule showed that a function can be viewed not only as a formula, but also as an ordered computation. Polynomial fitting then turned a prescribed representation into the matrix problem $Va\simeq y$, where least-squares approximation, interpolation and overparametrisation appear as different regimes of the same problem. Losses and gradients recast this fitting problem as optimisation. Training, mini-batches and epochs became names for familiar numerical operations, while the backward recursion through Horner's algorithm gave a first prototype of backpropagation. Neural networks then led from the fixed representation $V$ to a parametrised representation $V(\theta)$: the space in which the final approximation is performed is no longer prescribed in advance, but becomes part of what is learned.

This is the main change. In classical approximation, one often starts from a space and asks for the best element in that space. In learning, the data may also help construct the representation in which approximation becomes effective. The representation is no longer only the background of the problem; it becomes one of its unknowns. The approximating function is no longer only a member of a chosen space; it is the output of a parametrised computational graph.

The numerical questions therefore do not disappear. How is the function
represented? How is the error measured? Is the computation stable? Which
solution is selected when many solutions fit the data? How does the algorithm
influence the result? These questions become more challenging precisely because
representation, approximation, optimisation and selection are now coupled.

The same viewpoint suggests a way of formulating research questions. If a neural network is a trainable representation, then training does not merely select an element of a fixed approximation space; it changes the space in which approximation is performed. One may therefore ask convergence and stability questions not only for the final functions $F_{\theta,a}$, but also for the sequence of representations, or empirical spaces, generated during training. The passage from $V$ to $V(\theta)$ points toward a numerical analysis of learned approximation spaces: their approximation power, conditioning, stability, optimisation and selection. This connects naturally with nonlinear approximation \cite{DeVore1998}, with the additional feature that the approximating space is itself produced by a trainable computational graph.

This perspective is not restricted to fully connected feed-forward networks.
Different neural architectures may be interpreted as different rules for
constructing parametrised computational graphs
\cite{GoodfellowBengioCourville2016,HighamHigham2019}. The activation function
already changes the numerical character of the representation: with rectified
linear units, for instance, the computed functions are piecewise affine
\cite{GlorotBordesBengio2011}. Convolutional networks impose locality and
weight sharing through structured linear maps
\cite{LeCunBottouBengioHaffner1998}; residual networks alter the graph by adding
skip connections \cite{HeZhangRenSun2016}; recurrent networks use repeated maps
with shared parameters along a sequence \cite{Elman1990}; and graph neural
networks attach states to the nodes of a graph and update them through exchanges
along its edges \cite{ScarselliGoriTsoiHagenbuchnerMonfardini2009}. Transformers
give another important example: through attention, they produce state-dependent
nonlocal interactions \cite{VaswaniEtAl2017}. In this sense, an architecture is
a mathematical rule for generating a family of trainable representations.

\subsection{A teaching consequence}
This brings us back to the pedagogical motivation of the note. The aim is not to confine machine learning within numerical analysis. Probability and statistics remain central to questions of uncertainty, sampling and generalisation. The numerical viewpoint offers a complementary route. It starts from objects that many students have already met--polynomials, matrices, linear systems and derivatives--and organises them around representation, error, stability, optimisation and selection. In this way, machine learning can be introduced before the vocabulary and software become obstacles.

The regression line is a useful starting point, but it should not define the whole meaning of least squares. It is only the smallest example of approximation in a prescribed finite-dimensional space. A parabola, a cubic, a trigonometric polynomial or a spline leads to the same matrix structure and changes the representation without changing the basic problem. Least-squares approximation, interpolation and overparametrisation are then different regimes of one representational problem.

The terminology of learning can then be introduced naturally. A training set is a data set. Features are columns of an empirical representation. Parameters are unknowns to be selected. A loss is a residual measure. An optimiser is the numerical method used for that selection. If the computation is iterative, an epoch is a pass through the data, a batch is a subset of rows, and a learning rate is a step size. If the model is written as a computational graph and the optimiser requires gradients, backpropagation is the reverse application of the chain rule to that graph.

Horner's rule has a complementary role. It shows, in the simplest possible setting, that a function may be an ordered computation. In one variable this computation is a chain; in several variables it can be organised as a layered architecture of intermediate quantities. Nodes, layers, forward propagation and differentiation through a computation can therefore be introduced before neural networks are mentioned. The analogy has limits, but it gives the terminology a concrete mathematical origin.

This route is not limited to numerical analysis courses. Precisely because it starts from elementary mathematical objects, it can also serve courses in which mathematics is a tool rather than the main subject and, in simplified form, even the last years of secondary education. In this form, machine learning appears not as a departure from numerical mathematics, but as one of its contemporary continuations.

\subsection*{Acknowledgments}
The author acknowledges support from the Royal Society under the Wolfson Fellowship ``Uncertainty quantification, data-driven simulations and learning of multiscale complex systems governed by PDEs'' and from the FIS2023-01334 Advanced Grant ``Tackling complexity: advanced numerical approaches for multiscale systems with uncertainties'' (ADAMUS).

\end{document}